\documentclass[12pt]{amsart}
\usepackage{amssymb}
\usepackage{graphicx} 
\usepackage{epstopdf}

\newtheorem{thm}{Theorem}[section]
\newtheorem{lem}[thm]{Lemma}
\newtheorem{prop}[thm]{Proposition}
\newtheorem{cor}[thm]{Corollary}

\newtheorem{ddef}[thm]{Definition}

\def\De{\Delta}
\def\Ga{\Gamma}

\def\Z{{\mathbb{Z}}}
\def\N{{\mathbb{N}}}

\def\no={\,{\,|\!\!\!\!\!=\,\,}}

\def\({\left(}
\def\){\right)}

\def\lk{\operatorname{lk}}

\def\k{\mathbf k}
\def\ddim{\mathfrak{dim}}
\def\gin{\mathrm{gin}}
\def\l{\ell}

\def\beq{\begin{eqnarray}}
\def\eeq{\end{eqnarray}}

\begin{document}
\date{\today}

\title[Sequential Cohen-Macaulayness]{On Sequentially Cohen-Macaulay 
Complexes and
Posets}

\author[Bj\"orner]{Anders Bj\"orner}
\email{bjorner@math.kth.se}

\address{Department of Mathematics\\
          Royal Institute of Technology\\
          S-100~44 Stockholm, Sweden}

\author[Wachs]{Michelle Wachs$^1$}
\email{wachs@math.miami.edu}

\address{Department of Mathematics\\ University of Miami \\  Coral 
Gables, FL 33124,
USA}

\author[ Welker]{Volkmar Welker$^2$}
\email{welker@mathematik.uni-marburg.de}

\address{Fachbereich Mathematik und Informatik \\
          Universit\"at Marburg \\
          D-350~32 Marburg, Germany}

\thanks{1. Supported in part by National Science Foundation grants 
DMS DMS 0302310 and DMS 0604562.}
\thanks{2. Supported by Deutsche Forschungsgemeinschaft (DFG)}

\maketitle

\begin{abstract} The classes of sequentially Cohen-Macaulay and 
sequentially homotopy
Cohen-Macaulay complexes and posets are studied. First, some 
different versions of
the definitions are discussed and the homotopy type is determined. 
Second, it is shown how
various constructions, such as  join, product and rank-selection preserve these
properties. Third, a characterization of sequential 
Cohen-Macaulayness
for posets is given. Finally, in 
an appendix we outline
 connections with ring-theory and survey some 
uses of sequential Cohen-Macaulayness
 in commutative algebra.
\end{abstract}

\section{Introduction} \label{intro}

The notion of {\em sequential Cohen-Macaulayness}
is a nonpure generalization, due to Stanley \cite[Sec. III.2]{st}, of 
the notion of
{\em Cohen-Macaulayness}.  Stanley introduced this  in order to 
provide a  
ring-theoretic complement to the theory 
of nonpure shellability 
\cite{bw96,bw97}.   
Just as pure shellability implies 
Cohen-Macaulay\-ness,
nonpure shellability  implies sequential Cohen-Macaulay\-ness.

In this paper we show that the most common Cohen-Macaulay
preserving constructions on simplicial complexes and posets also preserve
sequential Cohen-Macaulayness.  This complements earlier results 
\cite{bw96,bw97} showing  that these
operations  preserve nonpure shellability.

We also discuss a nonpure version of Quillen's concept 
\cite{q} of 
{\em homotopy Cohen-Macaulayness},
introduced by the authors in \cite{bww}. This new concept is 
intermediate 
between nonpure shellability and
sequential Cohen-Macaulayness.
We show  that the homotopy version of sequential Cohen-Macaulayness
has the same strong topological consequences as that of (nonpure) 
shellability, namely
  having the  homotopy type of a wedge of spheres of (possibly) varying
dimensions.

Our primary goal is to extend results on Cohen-Macaulay simplicial complexes,
that have proven to be useful, to sequentially Cohen-Macaulay 
simplicial complexes.
Some of the proofs are straightforward
generalizations of the Cohen-Macaulay versions, while others 
require
substantially new ideas.
For basic facts from topological 
combinatorics we refer to the survey article by 
Bj\"orner \cite{bj2} 
and for background in commutative algebra to the books by 
Stanley 
\cite{st} and Bruns \& Herzog \cite{bh}.

 \smallskip
 
  As was mentioned, the notion of sequential Cohen-Macaulayness was 
 first defined in terms of
 commutative algebra by Stanley. In \cite{st} he also gave a homological characterization,
 see Appendix II, where this is outlined. Starting from Stanley's
 homological characterization, two other homological characterizations
 were found  by Duval \cite{D} and Wachs \cite{w1}. We take
 Wachs' characterization as our definition, and we will return
 to Duval's in the next section (see Proposition \ref{duv}).

Let $\Delta$ be a simplicial complex, and for
$0 \le m \le \dim \Delta$, let $\Delta^{\langle m\rangle}$ be the
subcomplex of $\Delta$ generated by its facets of dimension $\geq 
m$. 
 
 \begin{ddef}\label{def1} {\rm
 \begin{itemize}
 
\item[]
\item[(i)]
The complex
$\Delta$ is {\it sequentially acyclic} over $\k$ if $\Delta^{\langle 
m \rangle}$ is
$(m-1)$-acyclic over $\k$  for all $m = 0,1,\dots ,\dim \Delta$, i.e.,
$\tilde H_r(\Delta^{\langle m \rangle}; \k) = 0$ for all $r < m \le \dim
\Delta$, where $\k$ is the ring of integers or a field.
 \item[(ii)]
The complex $\Delta$ is {\it sequentially
connected} if
$\Delta^{\langle m \rangle}$ is
$(m-1)$-connected  for all $m = 0,1,\dots ,\dim \Delta$, i.e.,
$\pi_r(\Delta^{\langle m \rangle}) = 0$ for all $r < m \le \dim
\Delta$.  
\end{itemize} }
\end{ddef}
\noindent
  Recall that the {\it link} of a face $F$ in
$\Delta$ is defined to be the subcomplex
$$\lk_{\Delta}F = \{ G \in \Delta \mid F \cup G \in \Delta, \,\, F
\cap G = \emptyset\}. $$

  \begin{ddef}\label{def2} 
 \begin{itemize}
 \item[]
 \item[(i)]
The complex $\Delta$ is {\it sequentially Cohen-Macaulay}
  over $\k$ if $\lk_\Delta F$ is sequentially acyclic over $\k$ for all $F \in
\Delta$. \\ (Usually we will drop the reference to $\k$ and
just say ``sequentially acyclic'' and ``sequentially Cohen-Macaulay'' 
(SCM)). 
 \item[(ii)]
$\Delta$ is {\it sequentially homotopy Cohen-Macaulay} (SHCM)
if $\lk_\Delta F$ is sequentially connected for all $F \in \Delta$.
\end{itemize} 
\end{ddef}

A simplicial complex is said to be {\em pure} if all its facets
are of equal dimension. Clearly, a pure $d$-dimensional simplicial complex is
sequentially connected if and only if it is $(d-1)$-connected, and it is
sequentially acyclic if and only if it is
$(d-1)$-acyclic.  It follows that for pure simplicial complexes, the 
notion
 of ``S(H)CM'' reduces to
the notion of  ``(homotopy) Cohen-Macaulay''.

   A poset $P$ is said to be pure, sequentially connected, 
sequentially acyclic, SCM,
or SHCM if its order complex $\Delta(P)$ is, where $\Delta(P)$ is the 
simplicial
complex of chains of $P$.
\smallskip

The paper is organized as follows.  In Section~\ref{seqconn}, we show 
that a SHCM
complex has the homotopy type of a wedge of spheres in the dimensions 
of the facets
of the complex.  The homology version of this result appears in \cite{w1}.

The main result of
Section~\ref{constr} is that  sequential connectivity, sequential 
acyclicity, SHCMness,
and SCMness are all preserved by taking joins.   The proof of this   in 
the pure case
is quite simple, but in the nonpure case  relies on  a fiber lemma of
Quillen.  There are several interesting poset consequences of this result, in
particular a nice characterization of S(H)CM posets, and the result 
that  S(H)CMness is
preserved by  taking products of  posets with minimum elements.

In Section~\ref{rank} we prove some general results on induced subcomplexes of
S(H)CM complexes that enable us to show that  rank-selection on
semipure posets preserves the S(H)CM property and that  truncation on
general posets also preserves the S(H)CM property. 
 
Section~\ref{rankchar} contains a poset analog of Duval's characterization of
SCMness, which does not  follow directly from Duval's simplicial complex characterization.
This leads to a characterization of SCMness of
 semipure posets in terms of rank-selection, which extends a result
in the pure case due to 
Baclawski and Garsia \cite{BG} and Walker \cite{wal1}.  

Since Walker's  rank-selection result is unpublished and one of our results relies on this,
 we present his  proof
in an appendix.   In another appendix 
we outline the connection with 
commutative algebra. We
 give Stanley's definition of
sequentially Cohen-Macaulay modules and sketch the connection to the concept of
sequentially Cohen-Macaulay simplicial complexes as defined in this 
paper.
In this appendix we also survey some of the uses that the 
concept of
sequential Cohen-Macaulayness has found in the recent 
research literature.

\section{Sequentially connected complexes}\label{seqconn}

We begin with a description of the homotopy type of a sequentially connected
complex. The corresponding homological fact for sequentially acyclic 
complexes is
known from \cite{w1}.

\begin{prop}[{\cite[Proposition~1.9]{w1}}]  \label {sequential} Let
$\Delta$ be sequentially acyclic over {\bf k}.  Then $\tilde 
H_{\ast}(\De;\k)$ is free
 and
$\tilde H_{i}(\De ; \k)=0$ if $\Delta$ has no facet of dimension $i$.
\end{prop}

\begin{thm} Suppose that $\Delta$ is a simplicial complex with facets 
of dimensions
$d_1, \dots, d_t$.
  If $\Delta$ is sequentially connected, then it has the homotopy type 
of a wedge of
spheres of dimensions in $\{d_1, \dots, d_t\}$.
\end{thm}

\begin{proof} Assume that $d_1 > \dots > d_t$. We start with the case 
that $d_t \ge
2$, so that $\De =\Delta^{\langle d_t\rangle}$ is simply-connected. 
We already know
from  Proposition~\ref{sequential} that $\tilde H_{\ast}(\De;\Z)$ and $\tilde
H_{\ast}(\De^{\langle i \rangle};\Z)$ are free for all $i$. Let $\beta_i 
={\rm rank}\,
\tilde H_{i}(\De;\Z)= {\rm rank}\, \tilde H_{i}(\Delta^{\langle i\rangle};\Z)$.
We also know from  Proposition~\ref{sequential} that $\beta_i = 0$ 
for all $i < 2$.

   Since $\Delta^{\langle i\rangle}$ is $(i-1)$-connected,
  the Hurewicz theorem \cite[p. 479]{br} gives the existence of an isomorphism
$h_i: \pi_i(\Delta^{\langle i\rangle}) \to \tilde 
H_{i}(\Delta^{\langle i\rangle};\Z)$
when $i \ge 2$. This means that we can find mappings $\varphi_j^i$ \,
($j=1,\dots,\beta_i$) from the
$i$-sphere to $\Delta^{\langle i\rangle}$ whose induced homology 
classes form a basis
for the free group $\tilde H_{i}(\Delta^{\langle i\rangle};\Z)$. Let $W$ 
be a wedge of
spheres having $\beta_i$ \, $i$-dimensional spheres for all $i$. Since $\tilde
H_{i}(\De;\Z)= \tilde H_{i}(\Delta^{\langle i\rangle};\Z)$ we can piece the mappings
$\varphi_j^i$ together to a single mapping
$\Phi: W \to \De$ which induces isomorphism of homology in all 
dimensions. Since $W$
and $\De$ are simply-connected, the Whitehead theorem \cite[p. 
486]{br} implies that
such a mapping is a homotopy equivalence $\Phi: W\simeq \De$.

Assume now that $d_t=1$, and write $\De =\Delta^{\langle 2\rangle}
\cup\Ga$, where $\Ga$ is the $1$-skeleton of $\De$ (which is assumed to be
connected). Note that $\Delta^{\langle 2\rangle}$ (assumed to be 
simply-connected) is
by the preceding homotopy equivalent to a wedge of spheres. Let $T'$ be a
spanning-tree of the $1$-skeleton of
$\Delta^{\langle 2\rangle}$, and extend $T'$ to a spanning tree $T$ of $\Ga$.
Collapsing the tree $T$ turns $\De$ into a wedge of the space $\Delta^{\langle
2\rangle} / T'$ with a collection of loops ($1$-spheres), one coming 
from each edge
in $\Ga\setminus T$. Now, collapsing a contractible subspace does not 
change homotopy
type \cite[p. 436]{br}, so
$$ \De \simeq \De /T \simeq (\Delta^{\langle 2\rangle} /T')\,
\bigvee \{{\mathrm{loops}\}} \simeq \Delta^{\langle 2\rangle}
\bigvee \{{\mathrm{loops}}\}.$$

In the $d_t=0$ case there are in addition to the previous situation only some
isolated vertices. To handle these requires only adding the 
corresponding number of
$0$-spheres to the wedge already constructed.
\end{proof}

It was observed by Stanley \cite[p.~87]{st} that, just as in the pure 
case, SCMness is a
topological property, i.e., a property that depends only on the  geometric
realization of the simplicial complex and  ${\bf k}$; see 
\cite[Theorem
4.1.6]{w3}.  Although this is not true for the homotopy version (see 
\cite[Section
8]{q}), it is  easily seen to be true for both sequential 
connectivity and sequential
acyclicity.

\begin{prop} \label{topol} Sequential connectivity and sequential 
acyclicity are
topological properties.
\end{prop}

\begin{proof}  Given a nonnegative integer $m$ and a
topological space $X$, define
$X^{\langle m \rangle }$ to be the topological closure of the set
$$\{ p \in X : p \mbox{ has a
neighborhood  homeomorphic to}$$ \vspace{-.3in} $$\mbox{an open 
$d$-ball where $d \ge
m$}\}.$$     The result follows from the fact that
$\|\Delta^{\langle m\rangle}\| =
\|\Delta\|^{\langle m\rangle}$, for all $m = 0,1,\dots, \dim \Delta$, where
$\|\Delta\|$ denotes the geometric
realization of  $\Delta$.
\end{proof}

There is a characterization of SCMness due to Duval \cite{D} which 
involves the pure
$r$-skeleton of a simplicial complex. The {\it pure $r$-skeleton} 
$\Delta^{[r]}$ of a
simplicial complex
$\Delta$ is defined to be the subcomplex of $\Delta$ generated by  all faces of
dimension $r$.
In Proposition~\ref{duv} below we give  Duval's
formulation and a homotopy version of it.  It is shown in 
\cite[Theorem 1.5]{w1} that
Duval's formulation is equivalent to the one used here
(Definition 
\ref{def2}), 
by observing the
following connection between homology of the pure $r$-skeleton and sequential
acyclicity.  Recall that a simplicial complex
$\Delta$ is said to be {\it spherical} if it is $(\dim 
(\Delta)-1)$-connected. We say
that
$\Delta$ is  {\it homology-spherical} if it is
$(\dim (\Delta)-1)$-acyclic.

\begin{lem} \label{pureskel} A simplicial complex $\Delta$ is 
sequentially connected
(acyclic) if and only if its pure
$r$-skeleton
$\Delta^{[r]}$ is (homology-)spherical for all
$r
\le
\dim
\Delta$.
\end{lem}

\begin{proof} Let $\Delta^r$ denote the $r$-skeleton of $\Delta$.  Since
$\Delta^{[r]} = (\Delta^{<r>})^{r}$, we have $\pi_i(\Delta^{[r]} ) =
\pi_i(\Delta^{<r>})$ for all $i <r$.  Hence
$\Delta^{[r]}$ is $(r-1)$-connected if and only if
$\Delta^{<r>}$ is.  It follows that $\Delta$ is sequentially 
connected if and only if
all the pure skeleta of $\Delta$ are spherical.

An analogous argument works for homology.  \end{proof}

The homology version of the following result appears in \cite[Theorem 1.5]{w1}.

\begin{prop} \label{duv} A simplicial complex $\Delta$ is S(H)CM if 
and only if its
pure
$r$-skeleton $\Delta^{[r]}$ is (H)CM for all $r \le \dim \Delta$.
\end{prop}

\begin{proof} It is easy to see that if $F \in \Delta^{[ r]}$ then
$$\lk_{\Delta} (F)^{[r-\dim F -1]} =
\lk_{\Delta^{[ r] }}(F).$$ It therefore follows from Lemma~\ref{pureskel} that
$\lk_{\Delta^{[ r] }}(F)$ is spherical for all $r$ such that $F \in 
\Delta^{[ r]}$ if
and only if $\lk_{\Delta} (F) $ is sequentially connected.  This means that
$\Delta^{[ r] }$ is HCM for all $r$ if and only if $\lk_{\Delta} (F) $ is
sequentially connected for all $F$.

The analogous argument  works for homology  \cite{w1}.
\end{proof}

For more about the definition of sequential Cohen-Macaulayness, see
Appendix II.

\section{Join and product}\label{constr}

This section and the next deal with constructions on  complexes and posets
  that preserve SHCMness,  SCMness,
sequential connectivity, and sequential acyclicity.

Let  $P$ be a finite poset, and for
$x
\in P$ let $P_{\le x}$ denote the principal order ideal $\{y \in P \mid y \le
x\}$ and let $P_{< x} $ denote $\{y \in P \mid y <
x\}$.  For $x \le y \in P$, let $[x,y]$ denote the closed interval 
$\{z\in P \mid x
\le z\le y\}$ and let $(x,y)$ denote the open interval $\{z\in P \mid x
< z< y\}$.

The following  ``Quillen fiber lemma'' is the main tool in all of the 
results to
follow.

\begin{lem}\cite[Proposition 7.6]{q}. \label{quill}
  Let $f:P \to Q$ be a poset map,  and suppose that  the fiber
$\Delta(f^{-1}(Q_{\le q}))$  is $t$-connected for all $q \in Q$. 
Then $\Delta(P)$ is
$t$-connected if and only if $\Delta(Q)$ is $t$-connected.

The same is true with ``\,$t$-connected'' everywhere replaced by 
``\,$t$-acyclic''.
   \end{lem}

The following result is in the pure case an immediate consequence of 
the fact that the
join of an
$r$-connected complex and an $s$-connected complex is 
$(r+s-2)$-connected.   Its
proof in the nonpure case  relies on Lemma~\ref{quill}.  It is used 
to prove the
remaining results of this section.

\begin{thm} \label{join} The join of two sequentially connected 
(acyclic) simplicial
complexes is sequentially connected (acyclic).
\end{thm}

\begin{proof} Suppose $\Delta $ and
$\Gamma$ are sequentially connected.  We will show that
$(\Delta*\Gamma)^{\langle m \rangle}$ is $(m-1)$-connected for all $m$ by using
Lemma
\ref{quill}.
    Let $P$ be the poset of nonempty
faces of $(\Delta*\Gamma)^{\langle m \rangle}$ and let  $Q$ be the 
poset of nonempty
closed intervals of the totally ordered set $$(m-\dim \Gamma 
-1)<(m-\dim \Gamma ) <
\dots <\dim \Delta,$$ ordered by reverse inclusion. We construct a 
poset map from $P$
to
$Q$.    Let $F \in (\Delta*\Gamma)^{\langle m
\rangle}$.  Then
$F = F_1 \cup F_2$ where $F_1 \in \Delta$ and $F_2 \in \Gamma$.  Let $r =
\max \{\dim G \mid F_1 \subseteq G \in \Delta\}$ and $s = \max
\{\dim G \mid F_2 \subseteq G \in \Gamma\}$.  Define $f:P \to Q$ by $f(F) =
[m-r-1,s]$. It is easy to see that $f$ is order preserving, and  that 
for each $[a,b]
\in Q$, the fiber $f^{-1}(Q_{\le [a,b]})$ is the poset of nonempty faces of
$\Delta^{\langle m-a-1
\rangle} *
\Gamma^{\langle b \rangle}$.  It follows from the assumption that $\Delta$ and
$\Gamma$ are sequentially connected that  $\Delta^{\langle m-a-1 \rangle}$ is
$(m-a-2)$-connected and $\Gamma^{\langle b \rangle}$ is
$(b-1)$-connected.  Since the join of an
$r$-connected complex and an
$s$-connected complex is
$(r+s+2)$-connected, it follows that the order complex of the fiber 
$f^{-1}(Q_{\le
[a,b]})$ is
$(m-a-2+b-1+2)$-connected.  Since $m-a-2+b-1+2 \ge m-1$, the order 
complex of each
fiber is
$(m-1)$-connected.  Also
$\Delta(Q)$ is contractible since $Q$ has a minimum element.
Hence, Lemma \ref{quill} implies that
$\Delta(P)$ is
$(m-1)$-connected, which in turn implies that $(\Delta*\Gamma)^{\langle m
\rangle}$ is $(m-1)$-connected.

The homology version of the result is proved analogously.
\end{proof}

\begin{cor} The join of simplicial complexes $\Delta$ and $\Gamma$ is 
S(H)CM if and
only if $\Delta$ and $\Gamma$ are S(H)CM.
\end{cor}

\begin{proof} Use the fact that the link of a face in 
$\Delta*\Gamma$ is the join of
the links of faces in $\Delta$ and $\Gamma$.
\end{proof}

Recall that the {\em ordinal sum} of two posets $P$ and $Q$ is the poset on the disjoint union of $P$ and $Q$, whose order relation restricts to the ones on $P$ and $Q$ and sets all $p \in P$ below all $q \in Q$.
\begin{cor} \label{joinpos} The ordinal sum of posets is sequentially connected
(sequentially acyclic, SHCM, SCM) if and only if each poset is.
\end{cor}

   Since  the links of  faces of the order complex of a poset
$P$ are the order complexes of ordinal sums of   open intervals of
$\hat P$, where $\hat P$ is the poset $P$ with a minimum element $\hat 0$ and a
maximum element
$\hat 1$ attached,  we have the following  nice  characterization of 
S(H)CM posets.

\begin{cor} \label{int} A poset $P$ is S(H)CM  if and only if every open
interval of $\hat P$ is sequentially (connected)
acyclic.
\end{cor}

Given a simplicial complex $\Delta$, let  $P(\Delta)$ denote
the poset of nonempty faces of $\Delta$ ordered by inclusion.  Since
$\Delta(P(\Delta))$ is the barycentric subdivision of $\Delta$,  the 
two complexes
are homeomorphic.  We thus have the following
consequence of Proposition~\ref{topol} and Corollary~\ref{int}.

\begin{cor} A simplicial complex $\Delta$ is sequentially connected
(sequentially acyclic, SHCM, SCM) if and only if $P(\Delta)$ is 
sequentially connected
(sequentially acyclic, SHCM, SCM).
\end{cor}

\begin{proof} For sequential connectivity and acyclicity, this is an immediate
consequence of the fact that they are topological properties
(Proposition~\ref{topol}).  For SCMness it is also an immediate 
consequence of the fact
that SCMness is a topological property, which is more difficult to prove than
Proposition~\ref{topol}.  However, the result for both SCMness and SHCMness 
can be shown to
follow from the respective results for sequentially connectivity and sequential
acyclicity and from Corollary~\ref{int}.  Indeed, observe  that
$P(\lk_\Delta F)$ is isomorphic to the principal upper order ideal 
$\{G \in \Delta
\mid F
\subsetneq G\}$ of $P(\Delta)$.  Hence $\Delta$ is SHCM if and only 
if every open
principal upper order ideal of $P(\Delta)$ is sequentially connected. 
Each of the
other open intervals of
  $ P(\Delta)\cup \{\hat0,\hat 1\}$ is isomorphic to the proper part 
of a Boolean
algebra, and hence is sequentially connected no matter what $\Delta$ is.
Consequently
  $\Delta$ is SHCM if and only if every open interval of $ P(\Delta)\cup
\{\hat0,\hat 1\}$ is sequentially connected.  By Corollary~\ref{int}, 
$\Delta$ is SHCM
if and only if $P(\Delta)$ is.
\end{proof}

We  now turn to  the poset product operation. If $P_i$ is a finite 
poset, then let
$\acute P_i$ be the poset
$P_i$ with a minimum element $\hat 0_i$ attached,
and $\hat P_i$ be 
the poset $P_i$ with
a minimum element $\hat 0_i$ 
and a maximum element $\hat 1_i$ attached.

\begin{cor} \label{prodsc} Let $P_1$ and $P_2$ be sequentially 
(connected) acyclic
posets.  Then
\begin{enumerate}
\item[(1)]  $\acute P_1 \times \acute P_2 \setminus \{(\hat 0_1,\hat 0_2)\}$ is
sequentially (connected) acyclic,
\item[(2)] $\hat P_1 \times \hat P_2 \setminus\{(\hat 0_1,\hat 
0_2),(\hat 1_1,\hat
1_2) \}$ is sequentially (connected) acyclic.
\end{enumerate}
\end{cor}

\begin{proof} By a result of Quillen \cite{q} (see \cite[Theorem 5.1 
(b)]{wal2}),
there is a homeomorphism from
$\Delta( \acute P_1 \times \acute P_2 \setminus \{(\hat 0_1,\hat 0_2)\})$ onto
$\Delta(P_1) * \Delta(P_2)$.  Hence (1) follows from Theorems~\ref{topol}
and \ref{join}. Similarly, a homeomorphism of Walker
\cite[Theorem 5.1 (d)]{wal2} yields (2).
\end{proof}

\begin{cor} \label{prodcm} Let $P_1$ and $P_2$ be  posets with 
minimum elements
$\hat 0_1$ and $\hat 0_2$, respectively.  Then $P_1 \times P_2
$ is S(H)CM if and only if
$P_1$ and
$P_2$ are S(H)CM.
\end{cor}

\begin{proof} Suppose $P_1$ and $P_2$ are SHCM.  Then $P_i \setminus 
\{\hat 0_i\}$ is
SHCM.  By applying Corollary~\ref{prodsc} (1) to $P_i
\setminus \{\hat 0_i\}$ we get that $P_1 \times P_2 \setminus \{0_1,0_2\}$ is
sequentially connected, which implies that $P_1 \times P_2$ is sequentially
connected.  That all upper order ideals $(P_1 \times P_2)_{>(x_1,x_2)}$ are
sequentially connected also follows from from Corollary~\ref{prodsc} 
(1), and that all
open intervals
$((x_1,x_2),(y_1,y_2))$ of $P_1 \times P_2$ are sequentially 
connected follows from
Corollary~\ref{prodsc} (2). Now by Corollary~\ref{int}
$P_1 \times P_2$ is SHCM. Conversely, if $P_1 \times P_2$ is SHCM then the open
interval $(P_1
\times P_2)_{ >(x,\hat 0_2)}$, where $x$ is maximal in $P_1$, is 
SHCM.  Since this
open interval is isomorphic to $P_2 \setminus \{\hat 0_2\}$, we 
conclude that $P_2$
is SHCM.  Similarly $P_1$ is SHCM.

An analogous argument yields the homology version.
\end{proof}

The interval poset Int$(P)$ of a poset $P$ is the poset of closed 
intervals of $P$
ordered by inclusion.

\begin{cor} A poset $P$ is S(H)CM  if and only if \ {\rm Int}$(P)$ is S(H)CM.
\end{cor}

\begin{proof} The proof, which is similar to that of 
Corollary~\ref{prodcm}, uses the
homeomorphism in \cite[Theorem 6.1]{wal2}.
\end{proof}

\section{Rank selection}\label{rank}

Throughout this section, we assume that $P$ is a poset with
 a 
minimum element $\hat 0$.   The length $\ell(P)$ of $P$ is the length 
of
the longest chain of $P$, where   the length of a chain is
its cardinality minus one.
  For  $x \in P$, define the {\em rank},
$$r(x) := \ell([\hat 0, x]).$$  For $S \subseteq 
\{0,1,\dots,\ell(P)\}$, define the
rank-selected subposet
$$P_S := \{ x \in P \mid r(x) \in S\}.$$  
 
 It is well-known that
rank-selection preserves the Cohen-Macaulay property
 in the pure 
case, see e.g. \cite{BG} and
 \cite[p. 1858]{bj2} for references.  We 
shall show that the same is
true in the semipure case, where $P$ is said to be {\it semipure} if 
$[\hat 0, x]$ is
pure for all $x \in P$.  For  general bounded posets this is true 
only  for  special
types of rank-selection.

We begin with two lemmas.  The proof of the first lemma, given in
\cite{bj2}, relies on
Lemma~\ref{quill}.  It is used to prove the second lemma, which
is  key to the rest of this section.

  Given a simplicial complex
$\Delta$ and a subset $A$ of its vertex set, the induced subcomplex 
$\Delta(A)$ is
defined by
$$\Delta(A) := \{F \in \Delta \mid F \subseteq A\}.$$

\begin{lem}[{\cite[Lemma 11.11]{bj2}}] \label{induced} Let $\Delta$ 
be a simplicial
complex on vertex set $V$, and let $A \subseteq V$.  Assume that 
$\lk_\Delta(F)$ is
$m$-connected for all $F \in \Delta(V\setminus A)$.  Then $\Delta(A)$ is
$m$-connected if and only if $\Delta$ is  $m$-connected.

The same is true with ``\,$m$-connected'' everywhere replaced by 
``\,$m$-acyclic''.
\end{lem}

\begin{lem} \label{induced2} Let $\Delta$ be a S(H)CM  simplicial complex on
vertex set $V$.  Let  $A \subseteq V$ be such that for all facets
$F\in \Delta$,
$$|(V\setminus A)\cap F| = \begin{cases} 1  &\mbox{if } \dim F \ge t  \\ 0
&\mbox{otherwise}  \end{cases}$$ where $t$ is  some fixed element of 
$\{0,1,\dots
,\dim \Delta\}$. Then the induced subcomplex $\Delta( A)$ is S(H)CM.
\end{lem}

\begin{proof} We give the proof for the homotopy version.  An 
analogous argument
yields the homology version.  The proof is based on that of a pure 
version given in
\cite[Theorem~11.13]{bj2}.  

  We  begin by using Lemma~\ref{induced} to show that $\Delta(A)^{\langle
m\rangle}$ is $(m-1)$-connected for all $m=0,1,\dots,\dim \Delta(A)$. 
First note
that since $\Delta_{V\setminus A}$ is
$0$-dimensional, we only need to consider links of vertices that are 
not in $A$ when
applying Lemma~\ref{induced}.

{\bf Case 1.} $m < t$.  In this case \beq \label{eqinduc}\Delta(A)^{\langle
m\rangle} = \Delta^{\langle
m\rangle}(A).\eeq  Let $x \in V\setminus A$. Since no facet of 
dimension less than
$t$ contains $x$, we have
\begin{eqnarray*}\lk_{\Delta^{\langle
m\rangle}} (x) &=& \lk_{\Delta^{\langle
t\rangle}} (x)\\ &=& (\lk_\Delta x)^{\langle t-1\rangle}. \end{eqnarray*}
Since $\lk_\Delta x$ is sequentially connected, we have that 
$\lk_{\Delta^{\langle
m\rangle}} (x)$ is $(t-2)$-connected, which implies it is 
$(m-1)$-connected.   Since
$\Delta$ is sequentially connected, $\Delta^{\langle
m\rangle}$ is also $(m-1)$-connected.  Hence,  by Lemma~\ref{induced},
$\Delta^{\langle m\rangle}(A)$ is $(m-1)$-connected. By (\ref{eqinduc}),
$\label{eqind}\Delta(A)^{\langle m\rangle}$ is $(m-1)$-connected.

{\bf Case 2.} $m \ge t$.  In this case \beq \label{eqinduc2}\Delta(A)^{\langle
m\rangle} = \Delta^{\langle
m+1\rangle}(A).\eeq  Let $x \in V\setminus A$.  We have
$\lk_{\Delta^{\langle
m+1\rangle}} (x) = (\lk_\Delta x)^{\langle m\rangle}$, which is 
$(m-1)$-connected
since $\lk_\Delta x$ is  sequentially connected. Also $\Delta^{\langle
m+1\rangle}$ is $(m-1)$-connected for the same reason.  Hence by 
Lemma~\ref{induced},
$\Delta^{\langle
m+1\rangle}(A)$ is $(m-1)$-connected, which means by (\ref{eqinduc2}) that
$\Delta(A)^{\langle m\rangle}$ is $(m-1)$-connected.

We can now conclude that $\Delta(A)$ is sequentially connected.  To 
complete the
proof we need to show that all links of $\Delta(A)$ are sequentially connected.
This follows from the fact that for all $F \in \Delta(A)$,
$$\lk_{\Delta(A)} F = (\lk_\Delta F)(A).$$  Indeed, since the hypothesis of the
lemma holds for the sequentially Cohen-Macaulay complex $\lk_\Delta 
F$ and the set
$A$,   the  above argument yields the conclusion that $(\lk_\Delta F)(A)$ is
sequentially connected.
\end{proof}

A {\it completely balanced}  simplicial complex is defined to be a
simplicial complex $\Delta$ together with a
  ``coloring''
function  $\tau:V
\to \N$ such that for each facet
$F$, we have $\tau(F) =
\{0,1,2,\dots ,|F|-1\}$. The order complex of a semipure poset with  rank
function serving as color function is the prototypical
example a completely balanced simplicial complex.

\begin{thm} Let $\Delta$ be a completely balanced $d$-dimensional 
S(H)CM simplicial
complex.  Then for all $S \subseteq \{0,1,\dots,d\}$, the 
type-selected subcomplex
$\Delta_S:= \{G \in \Delta \mid \tau(G) \subseteq S\}$ is S(H)CM.
\end{thm}

\begin{proof} By induction, we can  assume that $|S| = d$.  Let
$t$ be the unique element of $ \{0,1,\dots,d\} -S$. Let $A$ be the 
set of vertices
whose color is not $t$.  Since $\Delta$ is completely balanced,  each facet of
dimension at least $t$ contains exactly one vertex of color $t$ 
(i.e., exactly one
vertex in the complement of $A$) and each facet of dimension less 
than $t$ contains
no vertices of color $t$.  The result now follows from Lemma~\ref{induced2}.
\end{proof}

In the case of posets, the
theorem reduces to the following:

\begin{cor} \label{rankselth} Let $P$ be a semipure S(H)CM poset of 
length $\ell$.
For all
$S \subseteq \{0,1,\dots,
\ell\}$, the rank-selected subposet $P_S$ is S(H)CM.
\end{cor}

\begin{thm} Let $P$ be a semipure S(H)CM poset of length $\ell$.  For all
$t= 0,1,\dots,\ell$,  the max-deleted subposet $$P^{(t)}:=P - \{x \in P \mid
x \mbox{ is a maximal element and }  r(x) \ge t \}$$ is S(H)CM.
\end{thm}

\begin{proof} Apply Lemma~\ref{induced2} to the order complex of $P$ 
by setting $A =
P^{(t)}$.
\end{proof}

\begin{cor} Let $\Delta$ be a S(H)CM simplicial complex. Then the
$t$-skeleton and $t$-coskeleton of $\Delta$ are S(H)CM for all
$t$.
\end{cor}

\begin{proof} The $t$-skeleton is obtained by rank-selection and the
$t$-coskeleton is obtained by repeated max-deletion in the face poset 
of $\Delta$.
\end{proof}

Now let us consider rank-selection in general posets.  Let $P$ be a 
bounded poset,
i.e., a poset with a minimum element and a maximum element.
  The {\it corank} $r^*(x)$ of $x\in P$
is  the rank of
$x$ in the dual poset $P^*$.  Note that when $P$ is pure,   $r^*(x) = 
\ell(P)-r(x)$.
For $S, T \subseteq \{0,1,\dots,\ell(P)\}$, define the {\it rank
selected subposet}
$$ P_S^T :=  \{x \in P \mid r(x) \in S \mbox{ and } r^*(x) \in T\}.$$

Rank selection in the general setting
does not preserve  S(H)CMness, see \cite[Figure~9a]{bw96}.  However,  a 
special type of
rank-selection, called truncation, does
preserve S(H)CMness, as can be seen in the next result.

\begin{thm} Let $P$ be a bounded S(H)CM poset of length $\ell$.  For 
$0 \le s,t <
\ell$, let $S = \{s,s+1,...,\ell\}$ and $T =\{t,t+1,\dots,\ell\}$.  Then the
truncation $P_S^T$ is S(H)CM.
\end{thm}

\begin{proof} Since $P_S^T = (((P_S)^*)_T)^*$,  we need only
prove the result for
$P_S$. Assume $s \ge 2$ (the result is trivial otherwise).  The poset 
$P_S$ can be
obtained by first removing
$\hat 0$ and then repeatedly removing all the atoms.  Apply 
Lemma~\ref{induced2} to
the order complex of
$P\setminus
\{\hat 0\}$ with
$A = (P\setminus
\{\hat 0\})
\setminus
\{ \mbox{atoms}\}$.
\end{proof}

\section{A poset analog of Duval's characterization}\label{rankchar}

  For a semipure poset $P$, let $P^{<j>}$ be the lower order ideal generated by all maximal elements
  of rank at least $j$ and let
  $P^{[j]}$ be the lower order 
ideal generated by
  all elements of rank $j$.  Clearly $\Delta(P^{<j>}) = \Delta(P)^{<j>}$.  Hence one has the following poset version of  Definition~\ref{def1}:  A semipure poset  $P$ is   is sequentially (connected) acyclic  if and only if $P^{<j>}$ is $(j-1)$-acyclic (-connected), for all $j=0,\dots, \l(P)$.  By combining this with Corollary~\ref{int}, one gets a nice characterization of S(H)CMness for semipure posets. 

  It is clearly {\em not} the case that $\Delta(P^{[j]}) = \Delta(P)^{[j]}$.  Hence a poset version of Duval's characterization of SCM (Proposition~\ref{duv}) does not follow directly from his characterization
  for simplicial complexes.  
Nevertheless a poset version of Duval's characterization does indeed hold.

\begin{thm}\label{char1}

A semipure poset $P$ is SCM (resp. SHCM) if and only if $P^{[j]}$ 
is
CM (resp. HCM) for all $j$.
\end{thm}

\begin{proof}
Suppose $P$ is S(H)CM.   Then $P^{<j>}$  is sequentially (connected) acyclic for all $ j$.  All intervals of $P^{<j>}$ are intervals of $P$ and all principal upper order ideals of $P^{<j>}$ have the form $I^{<k>}$,
where $I$ is a principal upper order ideal of $P$ and $k$ is some integer.  Hence by Corollary~\ref{int}, all intervals of $\widehat {P^{<j>}}$ are sequentially (connected) acyclic and therefore $P^{<j>}$ is S(H)CM.  It follows from Corollary~\ref{rankselth} that $P^{[j]} $ is (H)CM for all $j$ since $P^{[j]} = P^{<j>}_{\{0,1,\dots,j\}}$.

We now prove the homotopy version of the converse and leave the analogous proof of the  homology version to the reader.    Suppose $P$ semipure and $P^{[j]}$ is HCM for all $j$.  We will use Corollary~\ref{int} to show that $P$ is SHCM.  First we establish sequential connectivity for  $P$.   We  show that  $P^{<j>} $ is $(j-1)$-connected for all $j $ by induction on
$\l(P) - j$.   If $ j = \l(P)$  then $P^{<j>} = P^{[j]}$, which is $(j-1)$-connected.  Now assume  $j <\l(P)$.  If $P$ has no maximum element of rank $j$  then $P^{<j>} = P^{<j+1>}$, which by induction is $j$-connected and hence  $(j-1)$-connected.

Now assume $P$ has a maximal element of rank $j$.  Then $P^{<j>} = P^{<j+1>} \cup P^{[j]}$.  By induction,  $ P^{<j+1>} $ is $ j$-connected and hence $ (j-1)$-connected.  Also $P^{[j]}  $ is $(j-1)$-connected since it is HCM.  The intersection $P^{<j+1>} \cap P^{[j] }$ is the rank-selected subposet $(P^{[j+1]})_{\{0,1,\dots,j\} }$, which is HCM (by Theorem~\ref{rankselth}) and therefore $(j-1)$-connected.  By the Meyer-Vietoris 
Theorem~\cite[p. 229]{br},  $P^{<j>}$ is $(j-1)$-acyclic for all $j$, and by the Seifert - Van Kampen Theorem~\cite[p. 161]{br}, $P^{<j>}$ is simply connected.  It therefore follows from  the Hurewicz Theorem \cite[p. 479]{br} that $P^{<j>}$ is $(j-1)$-connected.  We can now conclude that $P$ is sequentially connected whenever $P$ is a semipure poset for which all the $P^{[j]}$ are HCM.

Next we  establish sequential connectivity for all intervals and principal upper order ideals of $P$.  Intervals of $P$ are intervals of some $P^{[j]}$, so they are indeed sequentially connected.   Let $I$ be a principal upper order ideal of $P$.  For each $j=0,\dots, \l(I)$,    the pure poset $I^{[j]}$  is a principal upper order ideal of $P^{[k]}$ where $k - j$ is the rank (in $P$) of the minimum element of $I$.   Since $P^{[k]}$ is HCM, so is the upper order ideal $ I^{[j]}$.  Since $I $ is a semipure poset for which all $I^{[j]}$ are HCM, by the result of  the previous paragraph, we conclude that  $I$ is sequentially connected.  It now follows from Corollary~\ref{int} that $P$ is SHCM.
\end{proof}

  It is known from work of Baclawski and Garsia \cite{BG} that 
Cohen-Macaulay
  posets are characterized by the property that all 
rank-selected subposets are
  homology-spherical. A sharpening of 
this result and of its homotopy version 
  was given by James W. 
Walker \cite{wal1}. The sharpening consists in the
  observation that 
it suffices to consider {\em intervals of ranks}.  
By an {\em interval of ranks} $S$ of a pure poset $P$ we mean a set 
of consecutive ranks
 $S=\{a, a+1, a+2, ...\} \subseteq \{1,2,\dots,\l(P \cup \{\hat 0\})\}$. 
If $P$ lacks a minimum element then by $P_S$ we mean $(P \cup \{\hat 0\})_S$.  (Recall rank-selection was defined in Section~\ref{rank} only for posets that have a minimum element.)

 \medskip

\begin{thm}[Walker \cite{wal1}] \label{walker1} 
If $P$ is a pure poset such that the rank-selected subposet $P_S$ 
is
spherical (resp. 
homology-spherical) for every
interval of ranks $S$, then $P$ is 
HCM (resp. CM).
\end{thm}

Since Walker's theorem is unpublished 
we give in Appendix I a version of his proof
which is distilled from his 
letter \cite{wal1}.

One might ask whether  Walker's theorem extends to the nonpure
setting, namely {\em if $P$ is semipure 
and $P_S$ is sequentially acyclic
for all sets of rank levels $S$, 
then is $P$ is SCM}. However, this is
false, as shown by the following 
counterexample.

\vspace*{.2in}
\begin{center}\includegraphics[width=1.5in]{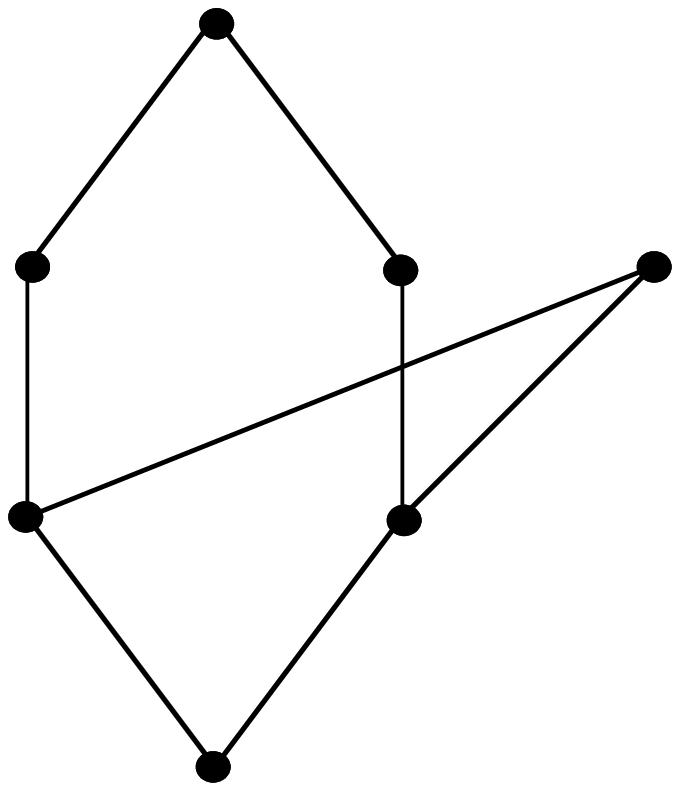}\end{center}

By combining
Walker's theorem and Theorem~\ref{char1} one gets a somewhat weaker extension of Walker's theorem, which is included in the following
list of characterizations of S(H)CMness.

\begin{cor}\label{char2}
Let $P$ be a semipure poset.  Then the following are equivalent.
\begin{enumerate}
\item $P$  is SCM (resp. SHCM). 
\item $(P^{[j]})_S$ 
is
$(|S|-1)$-acyclic (resp. $(|S|-1)$-connected) for all $j$ and 
 $S$ such that $S\subseteq [j]$.
 \item $(P^{[j]})_S$ 
is
$(|S|-1)$-acyclic (resp. $(|S|-1)$-connected) for all $j$ and 
every
interval of ranks $S$ such that $S\subseteq [j]$.
 \item $(P^{<j>})_S$ is $(|S|-1)$-acyclic (resp. $(|S|-1)$-connected) for all $j$ and 
 $S$ such that $S\subseteq [j]$.
\item $(P^{<j>})_S$ is $(|S|-1)$-acyclic (resp. $(|S|-1)$-connected) for all $j$ and 
 every
interval of ranks $S$ such that $S\subseteq [j]$.
\end{enumerate}
\end{cor}

\begin{proof} 
 $(1) \Rightarrow (2)$. Use Theorem~\ref{char1} and Corollary~\ref{rankselth}.

 $(2) \Rightarrow (3)$. Obvious.
 
 $(3) \Rightarrow (1)$.  Use Theorems~\ref{walker1} and~\ref{char1}.
 
 $(2) \iff (4)$ and $(3) \iff (5)$.  Use the fact that $(P^{[j]})_S= (P^{<j>})_S$ whenever $S\subseteq [j]$.
\end{proof}

\section*{Appendix I:  Walker's Proof}
\label{walkth}

\begin{proof}[Proof of Theorem~\ref{walker1}] We prove only the homotopy version, leaving the 
parallel
reasoning in the homology case to the reader.
 The argument 
proceeds in three stages, making use of the
following three 
properties of topological 
connectivity:
\begin{enumerate}
\item[(i)]If $X$ is an $n$-connected 
complex and $Y$ is an $(n-1)$-connected 
subcomplex, then $X/Y$ is 
$n$-connected.
\item[(ii)] Let $X= \bigvee X_i $ be a wedge of 
complexes. Then
$X$ is $n$-connected if and only if each $X_i$ is 
$n$-connected.
\item[(iii)] The suspension susp$(X)$ is $n$-connected 
if and only if
$X$ is $(n-1)$-connected.

\end{enumerate}

{\bf Step 
1:}\ {\em If $P$ is a pure poset such that $P_S$ is spherical for every lower interval of 
ranks $S$, and $a\in P$, then $P_{<a}$ is spherical.}

We 
may assume that $a$ is maximal in $P$. Let $n+1$ be the dimension of $X:=\Delta(P)$, and let $Y$ be the rank-selected subposet of $P$ formed by 
deleting the
level $P_{{\rm max}}$ containing $a$. Then $X$ is 
$n$-connected and $Y$ is
$(n-1)$-connnected, so $X/Y$ is 
$n$-connected. Since
$$ X/Y \cong \bigvee_{i\in P_{{\rm max}}}\ {\rm 
susp}(P_{< i})
$$
it follows that $  {\rm susp}(P_{< i})$ is 
$n$-connected for all $i\in P_{{\rm max}}$.
In particular, $P_{< a}$ 
is $(n-1)$-connected.

{\bf Step 2:}\ {\em If $P$ is a pure poset 
such that  $P_S$ is spherical for every  interval of 
ranks $S$, and $a\in P$, then $(P{<a})_D$ is spherical for
every upper interval $D$ of ranks of  $P_{<a}$.}

Suppose $D$ 
is an upper interval of ranks for $P_{<a}$, and let $\ell$ be
the 
index of the level of $P$ which contains $a$. Note that
$$(P_{<a})_D 
= (P_{D\cup \{\ell\}})_{<a}.
$$
The rank-selected subposet $P_{D\cup 
\{\ell\}}$ has the property that its rank-selected subposet  for
every interval of ranks is 
spherical, so $(P_{D\cup \{\ell\}})_{<a}$
is spherical by Step 
1.

{\bf Step 3:}\ {\em We complete the proof.}

If $b<a$ are 
elements of $P$, we want to show that the open interval $(b,a)$
is 
spherical. Observe that $(b,a)=(P_{<a})_{>b}$. By Step 2, $P_{<a}$ is 
a
ranked poset such that every upper rank-selected subposet is 
spherical.
Now apply the dual of Step 1.
\end{proof}

\section*{Appendix II: Sequential Cohen-Macaulayness in Commutative Algebra}
\label{commalg}

For a simplicial complex $\Delta$ over ground set $\{ 1, \ldots, n\}$, we
denote by $\k[\Delta]$ its Stanley-Reisner ring; that is the quotient 
of $S = \k[x_1, \ldots, x_n]$ by the ideal $I_\Delta$ generated by the
monomials $x_A = \prod_{i \in A} x_i$ in $S$ for $A \not\in \Delta$.
A simplicial complex $\Delta$ was defined by Stanley \cite{st} to be
sequentially Cohen-Macaulay over a field $\k$ if $\k[\Delta]$ is a sequentially
Cohen-Macaulay $S$-module. Again Stanley \cite[Definition 2.9]{st} defines a
{\em graded module} $M$ over a standard graded $\k$-algebra $R$ to be 
{\em sequentially
Cohen-Macaulay}, if there is a filtration of submodules
$$0 = M_0 \subset M_1 \subset \cdots \subset M_s = M$$
for which
\begin{itemize}
\item[$\triangleright$] $M_i/M_{i-1}$ is a Cohen-Macaulay $R$-module.
\item[$\triangleright$]  $\ddim\, M_1/M_0 < \cdots < \ddim\, 
M_s/M_{s-1}$, \\where $\ddim$ denotes Krull dimension.
\end{itemize}

The same concept appears in the work of Schenzel \cite{Sch}
under the name {\em Cohen-Macaulay filtered module}. Stanley 
mentions, and Schenzel
verifies \cite[Proposition 4.3]{Sch}, that if such a filtration exists then
it is unique. From \cite[Proposition 4.3]{Sch} it also follows that such a
filtration must coincide with the filtration defined as follows.
Let $N_j$ be the maximal submodule of $M$ of dimension $\leq j$, then
$N_{-1} = 0 \subseteq N_0 \subseteq \cdots \subseteq
  N_{\ddim\, M}=M$. If $M$ is sequentially Cohen-Macaulay
\cite[Proposition 4.3]{Sch} implies that for the indices
$-1 = j_{0} < \cdots < j_s = \ddim\, M$ for
which $N_{j_i-1} \neq  N_{j_i}$, $i \geq 1$,
we have $M_i = N_{j_{i}}$, $i \geq 0$.

In order to work out the maximal submodule of dimension $\leq j$ for a
Stanley-Reisner ring $\k[\Delta]$ we define $\Delta_j$ as
the subcomplex of $\Delta$ generated by its facets of dimension equal 
to $j$. 
(Recall that $\Delta^{\langle j \rangle}$ denotes  the
subcomplex of $\Delta$ generated by its facets of dimension $\geq 
j$). In general,
if $\Gamma$ is a subcomplex of $\Delta$ then we denote by
$I_{\Delta,\Gamma}$ the ideal generated in $\k[\Delta]$ by the
monomials $x_A$ with $A \in \Delta \setminus \Gamma$. As an $S$-module
$I_{\Delta,\Gamma}$ is isomorphic to $I_{\Gamma}/I_{\Delta}$.
Now if $\Delta$ is a simplicial complex and if $j_1 - 1 < \cdots < j_s -1$ are
the dimensions of the facets of $\Delta$ then it follows that
for $M = \k[\Delta]$ we have $M_{0} = 0$  and
$$M_{i-1} = I_{\Delta,\Delta^{\langle j_i-1\rangle}}, i \geq 1.$$
By \cite[III, Proposition 7.1]{st}
$$\ddim(I_{\Delta,\Delta^{\langle j_i-1 \rangle}}) = j_{i-1}.$$
One checks that $M_i$ is the maximal submodule of dimension $\leq j_i$ and
that if $N$ is a maximal submodule of dimension $\leq d$, for some $d$, then
$d = j_i$ for some $i$ and $N = M_i$.

Since
$$I_{\Delta,\Delta^{\langle j_i-1 \rangle}} \cong
   I_{\Delta^{\langle j_i-1 \rangle}} / I_\Delta$$ as $S$-modules, it 
follows that

$$M_{i-1}/M_{i-2} \cong
I_{\Delta^{\langle j_i-1 \rangle}}/I_{\Delta^{\langle j_{i-1}-1 
\rangle}} \cong $$

$$\cong
I_{\Delta_{j_{i-1}-1}, \Delta_{j_{i-1}-1}
\cap \Delta^{\langle j_i-1 \rangle}}.$$

The preceding isomorphism together with the definition of sequential
Cohen-Macaulayness and the fact that
$$\Delta^{\langle j \rangle} = \Delta_j \cup \cdots \cup 
\Delta_{\dim(\Delta)}$$
then yield the following characterization of simplicial
complexes $\Delta$ for which $\k[\Delta]$ is a sequential Cohen-Macaulay
$S$-module, given in \cite[III, Proposition 2.10]{st}:

{\em The $S$-module
$\k[\Delta]$ is sequentially Cohen-Macaulay, if and only if
$$I_{\Delta_i, \Delta_i \cap (\Delta_{i +1} \cup \cdots \cup 
\Delta_{\dim(\Delta)})}$$
is a Cohen-Macaulay module for all $i$.}

The latter condition translates into
$( \Delta_i , \Delta_i \cap
(\Delta_{i +1}  \cup \cdots \cup \Delta_{\dim(\Delta)}))$
being a relative simplicial complex which is Cohen-Macaulay over $\k$.
Stanley shows \cite[III, Theorem 7.2]{st} that a relative simplicial complex
$(\Delta, \Gamma)$ is Cohen-Macaulay over $\k$ if and only if for all
$A \in \Delta$ and all $i < \dim (\lk_\Delta(A))$ we have
$\widetilde{H}_i(\lk_\Delta(A), \lk_\Gamma(A);\k) = 0$.
This then gives rise to the characterizations of sequential Cohen-Macaulayness
by Duval \cite{D} (see Proposition \ref{duv}) and Wachs \cite{w1} (see Definition~\ref{def2} (i)).
Our definition of sequentially homotopy
Cohen-Macaulay  is a natural homotopy version of this formulation.
\smallskip

In commutative algebra the concept of a sequentially Cohen-Macaulay module
has been quite fruitful. We would like to mention
only a few of the developments which also bear a combinatorial flavor.
\begin{itemize}
\item
In \cite{hp} Herzog and Popescu characterize in ring-theoretic terms 
those 
 sequentially Cohen-Macaulay
$\k[\Delta]$ for which $\Delta$ is non-pure shellable.
This generalizes the concept of cleanness by Dress \cite{dr} which 
characterizes
the shellable $\Delta$ among the Cohen-Macaulay $\k[\Delta]$. Also in \cite{hp}
a proof of a homological characterization of sequentially 
Cohen-Macaulay modules,
originally due to Peskine \cite{st}, is given in terms of 
vanishing Ext-modules.
\item
It has
been shown in \cite{hrw} that $\k[\Delta]$ is sequentially 
Cohen-Macaulay if and only
if $I_{\Delta^*}$ has a componentwise linear resolution, where
$\Delta^* := \{ A ~|~\{1, \ldots, n\} \setminus A \not\in \Delta \}$ 
 is the ``combinatorial Alexander dual'' of $\Delta$.  This result
generalizes a result by Eagon and Reiner \cite{er} saying that $\k[\Delta]$ is
Cohen-Macaulay if and only if $I_{\Delta^*}$ has a linear 
resolution.
In \cite{hrw} it is also shown that these facts imply 
that duals
$\Delta^*$ of simplicial complexes $\Delta$ 
 for which $\k[\Delta]$ 
is sequentially Cohen-Macaulay
have the property that $\k[\Delta^*]$ is Golod. Despite substantial recent
progress (see \cite{bejo} and references therein) the latter property 
still waits for a
satisfying combinatorial characterization.
\item
In \cite{hs} it is shown that for an ideal $I$ in a polynomial ring
the Hilbert function of the local cohomology module of $R/I$ coincides
with the one of $R/ \gin(I)$ ($\gin(I)$ being the generic initial ideal
with respect to reverse lexicographic ordering) if and only if
$R/I$ is sequentially Cohen-Macaulay.
\end{itemize}

\end{document}